\documentclass[10pt]{amsart}
\usepackage{amsthm}
\usepackage{amsmath}
\usepackage{amssymb}
\usepackage{mathrsfs}
\usepackage{tikz-cd}
\usepackage{circledsteps}
\usepackage{caption}
\theoremstyle{plain}
\newtheorem{theorem}{Theorem}

\newtheorem*{covppthm}{Theorem}
\theoremstyle{definition}

\newtheorem{proposition}{Proposition}

\newtheorem{observation}[proposition]{Observation}
\newtheorem{corollary}[proposition]{Corollary}
\newtheorem{lemma}[proposition]{Lemma}
\newtheorem{definition}[proposition]{Definition}

\theoremstyle{remark}

\DeclareMathOperator{\prc}{prc}

\DeclareMathOperator{\tcf}{tcf}

\DeclareMathOperator{\norcf}{nor-cf}

\DeclareMathOperator{\bd}{bd}

\DeclareMathOperator{\pcfnor}{pcf_{nor}}

\DeclareMathOperator{\ran}{ran}

\DeclareMathOperator{\cf}{cf}

\DeclareMathOperator{\pcf}{pcf}
\renewcommand{\epsilon}{\varepsilon}

\DeclareMathOperator{\Reg}{{\sf Reg}}
\DeclareMathOperator{\NS}{NS}

\DeclareMathOperator{\cov}{cov}

\DeclareMathOperator{\pp}{pp}
\DeclareMathOperator{\ppnor}{pp_{nor}}

\newcommand{\sk}{\vskip.05in}

\numberwithin{equation}{section}
\numberwithin{proposition}{section}

\begin{document}
\title{Covering numbers and pcf theory}
\begin{abstract}
We develop a general technique for deriving covering theorems from pcf-theoretic information by refining sequences of partitions. The method gives a direct proof of Shelah’s $\operatorname{cov}$ vs. $\operatorname{pp}$ theorem and analogous results for diagonal coverings associated with normal ideals and coverings modulo arbitrary $\sigma$-complete ideals, where $\sigma$ is uncountable and regular. In each setting, we characterize the covering invariant through {pseudopowers,} true cofinalities of reduced products of progressive sets{,} and appropriate notions of cofinality in (possibly non-progressive) products of intervals of regular cardinals. 
\end{abstract}
\author{Todd Eisworth}
\date{\today}
\address{Department of Mathematics, Ohio University, Athens OH, 45701}
\thanks{This work is supported by the U.S. National Science Foundation under Award No.~2400200.}
\email{eisworth@ohio.edu}
\thispagestyle{empty}
\maketitle

\section{Introduction}

This paper develops an idea appearing as part of a proof from \cite{moreca} with the intent of providing a general technique for deriving bounds on covering numbers from pcf-theoretic information. The relationship between these two types of data is a major theme of pcf theory:  for example, the proof of Shelah's famous $\aleph_{\omega_4}$ theorem involves showing that the cofinality of the structure $([\aleph_\omega]^{\aleph_0},\subseteq)$ is equal to $\max\pcf\{\aleph_n:n<\omega\}$, with subsequent pcf analysis proving that this cardinal is less than $\aleph_{\omega_4}$. 

The covering problems under consideration all ask how large a family of small subsets of $\lambda$ must be in order to cover all subsets of a prescribed size in a certain way. Shelah has introduced notation for handling some of this complexity, and we recall it now:  

\begin{definition}
  Let $\sigma=\cf(\sigma)<\theta\leq\kappa\leq\lambda$ be infinite cardinals. A family $\mathscr{P}\subseteq[\lambda]^{<\kappa}$ is a \emph{$<\sigma$-cover of $[\lambda]^{<\theta}$} if every $u\in[\lambda]^{<\theta}$ is contained in the union of fewer than $\sigma$ members of $\mathscr{P}$. The minimum cardinality of such a family is denoted by $\cov(\lambda,\kappa,\theta,\sigma)$.
\end{definition}

As usual, for a set $X$ and a cardinal $\tau$, we write $[X]^{<\tau}$ and $[X]^{\leq\tau}$ for the subsets of $X$ of cardinality less than, or at most, $\tau$. The four parameters all play roles here, as we are asking how many sets from $[\lambda]^{<\kappa}$ it takes so that any member of $[\lambda]^{<\theta}$ is contained in some union of fewer than $\sigma$ of them.  There are also some familiar instances of this definition, for example,  $\cov(\lambda,\kappa^+,\kappa^+, 2)$ is just the cardinal $\cf([\lambda]^\kappa,\subseteq)$.  Another common setting is $\cov(\lambda,\kappa^+,\theta^+,\sigma)$ with $\theta\leq\kappa<\lambda$. With this configuration,  the  sets from our covering family are allowed to have cardinality at most $\kappa$, and the requirement is that every set of cardinality at most $\theta$ can be covered using fewer than $\sigma$ of them.

There is an easy reason to expect a connection between these covering numbers and products of regular cardinals. Let us define
$$\mathsf{R}(\kappa,\lambda]=(\kappa,\lambda]\cap\Reg$$
as our notation for the set of regular cardinals in the interval $(\kappa,\lambda]$. Suppose $X\in[\lambda]^{\leq\kappa}$ and $A\subseteq \mathsf{R}(\kappa,\lambda]$. We write $\prod A$ for the set of functions $f$ with domain $A$ such that $f(a)<a$ for every $a\in A$. Since $|X|\leq\kappa<a$ and $a$ is regular, the function
$$h_X(a)=\sup(X\cap a)+1$$
belongs to $\prod A$. Moreover, if $\gamma\in X\cap a$, then $\gamma<h_X(a)$. Covering a set of ordinals by a few small sets therefore gives bounds for those ordinals in the corresponding coordinates of a product. The main task is to go in the other direction and use bounds in products to construct coverings, and this paper exploits one natural way of doing this.

We recall the pcf notation needed to make this correspondence precise. Our conventions mostly follow~\cite{abrahammagidor} and~\cite{kojmanabc}. For more obscure material, we keep close to~\cite{cardarith} and note that the analytical guide~\cite{analyticalguide} gives a useful account of the relevant parts of the book and later developments. The surveys~\cite{burkemagidor,jechpcf} provide other introductions to pcf theory.

A set $A$ of regular cardinals is \emph{progressive} if
$$|A|<\min(A).$$
For a proper ideal $J$ on $A$, we write $J^+$ for the collection of $J$-positive subsets of $A$, and use the usual ordering modulo $J$:
$$f\leq_J g\Longleftrightarrow\{a\in A:f(a)>g(a)\}\in J.$$
Similarly, $f<_J g$ means that $\{a\in A:f(a)\geq g(a)\}\in J$. The ideal specifies which sets of coordinates may be ignored when comparing functions. If there is a regular cardinal $\tau$ and a $<_J$-increasing cofinal sequence
$$\langle f_\alpha:\alpha<\tau\rangle$$
in $\prod A/J$, then the \emph{true cofinality} of the product is $\tau$, written
$$\tau=\tcf(\prod A/J).$$
Such a sequence is called a \emph{scale modulo $J$}. 

For progressive $A$, $\pcf(A)$ is the set of possible cofinalities of ultraproducts of $A$. For a cardinal $\tau$, the ideal
$$J_{<\tau}[A]=\{B\subseteq A:\pcf(B)\subseteq\tau\}$$
consists of those subsets whose possible cofinalities are all less than $\tau$. The following standard facts about pcf generators are the ones we will use; see~\cite{cardarith} and Section~5 of~\cite{analyticalguide}.

\begin{proposition}
  Let $A$ be a nonempty progressive set of regular cardinals. Then $\pcf(A)$ has a largest member $\Lambda$, and there is a generating sequence
  $$\langle B_\tau[A]:\tau\in\pcf(A)\rangle$$
  with $B_\Lambda[A]=A$ such that, for every $\tau\in\pcf(A)$,
  $$J_{<\tau^+}[A]=J_{<\tau}[A]+B_\tau[A].$$
  Moreover, $\prod B_\tau[A]$ has a pointwise cofinal family of cardinality $\tau$.

  If $J$ is a proper ideal on $A$ with $J_{<\tau}[A]\subseteq J$ and $B_\tau[A]\notin J$, then
  $$\tau=\tcf\bigl(\prod B_\tau[A]/(J\upharpoonright B_\tau[A])\bigr).$$
  Conversely, if $\tau=\tcf(\prod A/J)$, then $\tau\in\pcf(A)$, $J_{<\tau}[A]\subseteq J$, and $B_\tau[A]\notin J$.
\end{proposition}

Here $J_{<\tau}[A]+B_\tau[A]$ is the ideal generated by $J_{<\tau}[A]$ together with $B_\tau[A]$. More generally, $J_{<\tau}[A]$ is generated by the sets $B_\rho[A]$ for $\rho\in\pcf(A)\cap\tau$. Consequently, given any proper ideal $J$ on $A$, we may choose the least $\tau$ for which $B_\tau[A]\notin J$ and obtain true cofinality $\tau$ by restricting to that generator. Taking $\tau=\Lambda$ also gives the familiar equality $\cf(\prod A,\leq)=\max\pcf(A)$. We fix a generating sequence for each progressive set under consideration.

For the covering problems we have in mind, we are going to need to consider several variations of the idea of cofinality in a product of cardinals. The following pins down one such notion of cofinality, in a context that will allow us to talk about Shelah's previous work. 

\begin{definition}
  Let $\sigma=\cf(\sigma)<\theta$ be infinite cardinals, and let $A$ be a set of regular cardinals. 
  \begin{enumerate}
    \item A family $\mathscr{F}\subseteq\prod A$ is \emph{$\Gamma(\theta,\sigma)$-cofinal in $\prod A$} if, whenever $B\in[A]^{<\theta}$ and $g\in\prod B$, there is a family $\mathscr{F}_0\in[\mathscr{F}]^{<\sigma}$ such that for every $b\in B$ there is an $f\in\mathscr{F}_0$ with $g(b)<f(b)$. The minimum cardinality of such a family is denoted by $\cf_{\Gamma(\theta,\sigma)}(\prod A)$.

  \sk 
\item If $A$ is progressive, then $\pcf_{\Gamma(\theta,\sigma)}(A)$ is the set of all cardinals of the form $\tcf(\prod B/J)$, where $B\in[A]^{<\theta}$, $J$ is a proper $\sigma$-complete ideal on $B$, and the true cofinality exists.
\end{enumerate}
\end{definition}

The family $\mathscr{F}_0$ may depend on both $B$ and $g$, and different members of $\mathscr{F}_0$ may provide the bound at different coordinates. The first parameter of $\Gamma(\theta,\sigma)$ controls the size of $B$, and the second controls the size of $\mathscr{F}_0$. The $\sigma$-completeness condition on the ideals gives the corresponding restriction on the quotient side: a union of fewer than $\sigma$ exceptional sets still belongs to the ideal. Lemma~\ref{conventional} shows that the minimum size of a $\Gamma(\theta,\sigma)$-cofinal family is exactly $\sup\pcf_{\Gamma(\theta,\sigma)}(A)$ when $A$ is progressive.

We can now express the elementary connection with coverings more precisely. If $\mathscr{P}$ witnesses $\cov(\lambda,\kappa^+,\theta^+,\sigma)$ and $A\subseteq \mathsf{R}(\kappa,\lambda]$, then
$\{h_X:X\in\mathscr{P}\}$
is $\Gamma(\theta^+,\sigma)$-cofinal in $\prod A$. To see this, take $B\in[A]^{\leq\theta}$ and $g\in\prod B$, and cover the range of $g$ by fewer than $\sigma$ members of $\mathscr{P}$. For every $b\in B$, one of these sets $X$ contains $g(b)$, and hence $g(b)<h_X(b)$.

To describe the same information in terms of singular cardinals, we use pseudopowers. These collect true cofinalities of products of regular cardinals approaching a fixed singular cardinal, with restrictions on the size of the product and the completeness of the ideal.

\begin{definition}
  Let $\nu$ be singular and suppose $\sigma=\cf(\sigma)\leq\cf(\nu)<\theta$. Then $\pp_{\Gamma(\theta,\sigma)}(\nu)$ is the supremum of all cardinals of the form $\tcf(\prod B/J)$, where $B\subseteq\nu\cap\Reg$ is cofinal in $\nu$, $|B|<\theta$, and $J$ is a proper $\sigma$-complete ideal on $B$ containing every bounded subset of $B$. We write $J^{\bd}[B]$ for the ideal of bounded subsets of $B$.
\end{definition}

The requirement $J^{\bd}[B]\subseteq J$ ensures that discarding a bounded initial segment of $B$ does not change the true cofinality.  The restriction $\sigma\leq\cf(\nu)$ is necessary: otherwise, $B$ would be the union of fewer than $\sigma$ bounded sets, and no proper $\sigma$-complete ideal on $B$ could contain all of them.

The state of the art for these covering numbers is Shelah's $\cov$ vs. $\pp$ theorem (Theorem~5.4 of Chapter~II of~\cite{cardarith}), which we state here in the configuration of parameters used in this paper.

\begin{covppthm}[Shelah's $\cov$ vs. $\pp$ theorem]
  Suppose $\aleph_0<\sigma\leq\theta\leq\kappa<\lambda$ with $\sigma$ regular. Then
  \begin{equation*}
    \cov(\lambda,\kappa^+,\theta^+,\sigma)
    =\lambda+\sup\{\pp_{\Gamma(\theta^+,\sigma)}(\nu):
    \sigma\leq\cf(\nu)\leq\theta,\ \kappa<\nu\leq\lambda\}.
  \end{equation*}
\end{covppthm}

Further results concerning the theorem appear in~\cite{furtherca}. Our prior work in the area gave a second proof of Shelah's $\cov$ vs. $\pp$ theorem in~\cite{eisworthcompactness}, and further developed the theory of representation to establish a combinatorial principle we called the Pseudopower Dichotomy in~\cite{eisworthdichotomy}. The work in this paper should be viewed as a continuation of this project.  We show that our techniques quickly imply the $\cov$ vs. $\pp$ theorem, and then we deploy them to obtain some new results.

Section~\ref{sec:partitions} develops the partition refinement machinery used throughout the paper, including the elementary closure observation. In Section~\ref{sec:coveringnumbers} we apply this machinery to give a new and direct proof of Shelah's $\cov$ vs. $\pp$ theorem, together with four characterizations of the covering number. Section~\ref{sec:normalideals} develops a corresponding theorem for normal ideals and diagonal coverings, and Section~\ref{sec:almosteverywhere} treats coverings modulo an arbitrary $\sigma$-complete ideal. In all of these settings, we will obtain the same four types of characterizations of the cardinal under consideration.

\section{Coverings from partitions}
\label{sec:partitions}

In this section we develop the general technique that underlies our results. The two main ingredients are Lemma~\ref{partitiontree}, which turns suitable trees of refinements into coverings, and Lemma~\ref{partitioncriterion}, which uses bounds in products of regular cardinals to supply the required refinements. We finish with an elementary closure construction that controls the size of the family of partition sequences.

\begin{definition}
  Let $\kappa$, $\lambda$, and $\theta$ be infinite cardinals.
  \begin{enumerate}
    \item A $(\kappa,\theta)$-partition sequence of $\lambda$ is a sequence $\bar{c}=\langle c_\beta:\beta<\theta\rangle$ of functions from $\lambda$ into $\kappa$.  Each $c_\beta$ should be thought of as a partition of $\lambda$ into $\kappa$ cells (some of which may be empty), and the full partition sequence can be visualized as a $\kappa\times\theta$ matrix of subsets $X_{\alpha\beta}$ of $\lambda$, where $\alpha<\kappa$, $\beta<\theta$, and
      \begin{equation}
        \label{matrixpiece}X_{\alpha\beta}=\{\gamma<\lambda: c_\beta(\gamma)=\alpha\}.
      \end{equation}

      \medskip

    \item If $\bar{c}$ and $\bar{d}$ are $(\kappa,\theta)$-partition sequences of $\lambda$, then $\bar{d}$ \emph{refines} $\bar{c}$ if, for every $\beta<\theta$, the partition in column $\beta$ of $\bar{d}$ refines the partition in column $\beta$ of $\bar{c}$.
  \end{enumerate}
\end{definition}

Fix a $(\kappa,\theta)$-partition sequence $\bar{c}=\langle c_\beta:\beta<\theta\rangle$ of $\lambda$ and a function $g:\theta\rightarrow\lambda$. For each $\beta<\theta$, let $g^{\bar{c}}(\beta)$ be the unique cell in column $\beta$ containing $g(\beta)$. Thus,
$$g^{\bar{c}}(\beta)=X_{c_\beta(g(\beta)),\beta}=\{\gamma<\lambda:c_\beta(\gamma)=c_\beta(g(\beta))\},$$
and
$$g\in\prod_{\beta<\theta}g^{\bar{c}}(\beta).$$
If $\bar{d}$ refines $\bar{c}$, then $g^{\bar{d}}(\beta)\subseteq g^{\bar{c}}(\beta)$ for every $\beta<\theta$. We will use cofinal families in products of regular cardinals to shrink cells whose cardinality exceeds $\kappa$. The next definition identifies the regular cardinals involved.

\begin{definition}
  Let $\bar{c}=\langle c_\beta:\beta<\theta\rangle$ be a $(\kappa,\theta)$-partition sequence for $\lambda$ with $\theta\leq\kappa<\lambda$ and $X_{\alpha\beta}$ defined as in (\ref{matrixpiece}).
  \begin{enumerate}
    \item For $\alpha<\kappa$ and $\beta<\theta$, we define
      $$\lambda^{\bar{c}}_{\alpha\beta}=|X_{\alpha\beta}|,$$
      and
      $$\kappa^{\bar{c}}_{\alpha\beta}=\cf(\lambda^{\bar{c}}_{\alpha\beta}).$$
      We may omit $\bar{c}$ if it is clear from context.

      \medskip

    \item These cardinals classify the indices in $\kappa\times\theta$ as \emph{small}, \emph{medium}, or \emph{large}:
      \begin{enumerate}
          \sk
        \item $(\alpha,\beta)$ is \emph{small} if $\lambda_{\alpha\beta}\leq\kappa$;
          \sk
        \item $(\alpha,\beta)$ is \emph{medium} if $\kappa_{\alpha\beta}\leq\kappa<\lambda_{\alpha\beta}$; and
          \sk
        \item $(\alpha,\beta)$ is \emph{large} if $\kappa<\kappa_{\alpha\beta}$.
          \sk
      \end{enumerate}
      When several partition sequences are under consideration, we may say, for example, that $(\alpha,\beta)$ is \emph{small in $\bar{c}$}. We also call the corresponding cell $X_{\alpha\beta}$ small, medium, or large.
  \end{enumerate}
\end{definition}

The idea is to start with a function $g:\theta\rightarrow\lambda$ and some $(\kappa,\theta)$-partition sequence for $\lambda$, and do successive refinements of a certain form. Our refining process treats the three classes of cells discussed above in different ways, and it is our pcf assumptions that allow us to deal sensibly with the large cells. We iterate our refinement process through countably many stages: since cell cardinalities cannot decrease indefinitely, suitable refinements eventually place the relevant values of $g$ in small cells.   

When we refine a partition sequence, we want to make sure that in each column the resulting partition has at most $\kappa$ cells, but we also want to make sure that we are breaking cells up into smaller pieces where possible.   The small cells already have the required size for the type of covering we seek and are left alone. The medium cells can be partitioned into at most $\kappa$ sets of smaller cardinality because of our cofinality assumptions, so we can break them apart accordingly.

For a large cell, we need to do something smarter as we cannot just break it into pieces of smaller cardinality: this would result in a partition with more than $\kappa$ cells. Instead, we will use a function in a product to tell us how we carve off a suitable smaller set.  To see how this works, note that the large cells determine a natural set of regular cardinals: for a partition sequence $\bar{c}$, define
$$A^{\bar{c}}=\{\kappa_{\alpha\beta}: (\alpha,\beta)\text{ large in }\bar{c}\}.$$
Then $A^{\bar{c}}\subseteq \mathsf{R}(\kappa,\lambda]$ and $|A^{\bar{c}}|\leq\kappa$, so $A^{\bar{c}}$ is a progressive set of regular cardinals whenever it is nonempty, and we can use our pcf assumptions applied to $\prod A^{\bar{c}}$. (Be warned that in the future we will omit the superscript $\bar{c}$ when the partition sequence is fixed.)  With this goal in mind, the following definition makes sense.

\begin{definition}
  Let $\bar{c}$ be a $(\kappa,\theta)$-partition sequence of $\lambda$, where $\theta\leq\kappa<\lambda$, and let $g:\theta\rightarrow\lambda$. A collection $\{\bar{d}_\zeta:\zeta<\zeta^*\}$ of $(\kappa,\theta)$-partition sequences \emph{reduces $g$ over $\bar{c}$} if the following conditions hold:
  \begin{itemize}
    \item each $\bar{d}_\zeta$ is a refinement of $\bar{c}$, and

      \sk

    \item for every $\beta<\theta$, if $g^{\bar{c}}(\beta)$ is not small in $\bar{c}$, then there is a $\zeta<\zeta^*$ such that $|g^{\bar{d}_\zeta}(\beta)|<|g^{\bar{c}}(\beta)|$.
  \end{itemize}
\end{definition}

Already we can see that this is useful! Suppose the singleton collection $\{\bar{c}_{n+1}\}$ reduces $g$ over $\bar{c}_n$ for every $n<\omega$. If $C_n$ is the union of the small cells of $\bar{c}_n$, then $\bigcup_{n<\omega}C_n$ covers the range of $g$. Otherwise, for some $\beta<\theta$, the sequence $\langle|g^{\bar{c}_n}(\beta)|:n<\omega\rangle$ would be strictly decreasing. Different refinements may shrink the cells containing different values of $g$, so the next lemma extends this argument to one organized along a tree.  We will use this lemma repeatedly: it supplies the mechanism for building covering families from partition sequences.

\begin{lemma}[Tree covering lemma]
  \label{partitiontree}
  Suppose $\theta\leq\kappa<\lambda$, and $T$ is a non-empty tree of finite sequences ordered by extension and closed under initial segments. Let $\langle\bar{c}^s:s\in T\rangle$ be a collection of $(\kappa,\theta)$-partition sequences and assume that $\bar{c}^t$ refines $\bar{c}^s$ whenever $t$ extends $s$ in $T$. Let $C^s$ be the union of the small cells of $\bar{c}^s$, let $g:\theta\rightarrow\lambda$, and fix $\beta<\theta$. Suppose that whenever $s\in T$ and $g^{\bar{c}^s}(\beta)$ is not small, there is an immediate successor $t\in T$ of $s$ such that
  $$|g^{\bar{c}^t}(\beta)|<|g^{\bar{c}^s}(\beta)|.$$
  Then $g^{\bar{c}^s}(\beta)$ is small for some $s\in T$, and in particular
  $g(\beta)\in\bigcup_{s\in T}C^s.$
\end{lemma}
\begin{proof}
  Suppose that $g^{\bar{c}^s}(\beta)$ is not small for any $s\in T$. By hypothesis, every $s\in T$ has an immediate successor $t\in T$ satisfying
  $$|g^{\bar{c}^t}(\beta)|<|g^{\bar{c}^s}(\beta)|.$$
  Starting from the root, we may recursively choose an infinite branch $\langle s_n:n<\omega\rangle$ through $T$ such that
  $$|g^{\bar{c}^{s_{n+1}}}(\beta)|<|g^{\bar{c}^{s_n}}(\beta)|$$
  for each $n<\omega$. This gives a strictly decreasing sequence of cardinals, a contradiction.
\end{proof}

Our next move is to develop a way of producing partition sequences that satisfy the assumptions of the preceding lemma.  To accomplish this goal,  we fix for each $X\subseteq\lambda$ of cardinality greater than $\kappa$ a filtration $\langle X[\epsilon]:\epsilon<\cf(|X|)\rangle$ such that
\begin{itemize}
  \item the sequence $\langle X[\epsilon]:\epsilon<\cf(|X|)\rangle$ is increasing and continuous with union $X$, and

    \sk

  \item $|X[\epsilon]|<|X|$ for each $\epsilon$.

    \sk
\end{itemize}
We do not need to make any special assumptions about these filtrations; all that matters is that we fix the assignment
$$X\longmapsto\langle X[\epsilon]:\epsilon<\cf(|X|)\rangle$$
as a parameter. Each filtration induces a partition of $X$ by assigning every $\delta\in X$ to the least $\epsilon$ such that $\delta\in X[\epsilon]$.

\begin{definition}
  Let $\bar{c}$ be a $(\kappa,\theta)$-partition sequence of $\lambda$, set $A=A^{\bar{c}}$, and let $h\in\prod A$. The \emph{$h$-refinement of $\bar{c}$}, denoted $\bar{c}\downarrow h$, is constructed as follows:
  \begin{itemize}

    \item If $(\alpha,\beta)$ is small in $\bar{c}$, then $X_{\alpha\beta}$ remains unchanged as a cell in column $\beta$ of $\bar{c}\downarrow h$.
      \sk

    \item If $(\alpha,\beta)$ is medium in $\bar{c}$, then $X_{\alpha\beta}$ is replaced by the at most $\kappa$ cells in the partition induced by its filtration.

      \sk

    \item If $(\alpha,\beta)$ is large in $\bar{c}$, then $X_{\alpha\beta}$ is split into the two sets
      $$X_{\alpha\beta}[h(\kappa_{\alpha\beta})]$$
      and
      $$X_{\alpha\beta}\setminus X_{\alpha\beta}[h(\kappa_{\alpha\beta})].$$
  \end{itemize}
\end{definition}
There are at most $\kappa$ resulting cells in each column. We fix reindexings by ordinals below $\kappa$, so that $\bar{c}\downarrow h$ is again a $(\kappa,\theta)$-partition sequence and refines $\bar{c}$. In particular,
$$g^{\bar{c}\downarrow h}(\beta)\subseteq g^{\bar{c}}(\beta)$$
for every $\beta<\theta$. Every medium cell is split into cells of smaller cardinality. To obtain a smaller cell containing $g(\beta)$ when $X_{\alpha\beta}$ is large, it suffices for $h(\kappa_{\alpha\beta})$ to exceed the least $\epsilon$ such that $g(\beta)\in X_{\alpha\beta}[\epsilon]$. The next definition packages these requirements into a single function.

\begin{definition}
  \label{partitionbound}
  Let $\bar{c}$ be a $(\kappa,\theta)$-partition sequence of $\lambda$, where $\theta\leq\kappa<\lambda$, let $A=A^{\bar{c}}$, and let $g:\theta\rightarrow\lambda$. Define $\hat{g}\in\prod A$ by letting $\hat{g}(\eta)$ be the least $\epsilon<\eta$ such that, for every $\alpha<\kappa$ and $\beta<\theta$, if $\eta=\kappa_{\alpha\beta}$ and $g(\beta)\in X_{\alpha\beta}$, then $g(\beta)\in X_{\alpha\beta}[\epsilon]$.
\end{definition}
The function $\hat{g}$ is well-defined. For a fixed $\eta\in A$, at most one cell in each column contains $g(\beta)$, so there are at most $\theta$ relevant filtration stages. Since $\eta$ is regular and $\theta\leq\kappa<\eta$, those stages are bounded below $\eta$. Moreover, $\hat{g}$ is zero except possibly on a set of at most $\theta$ cardinals of the form $\cf(|g^{\bar{c}}(\beta)|)$.  Our next lemma  shows us how to translate bounds on $\hat g$ into the shrinking refinements needed to apply Lemma~\ref{partitiontree}.

\begin{lemma}[Refinement criterion]
  \label{partitioncriterion}
  Suppose $\bar{c}$ is a $(\kappa,\theta)$-partition sequence of $\lambda$, where $\theta\leq\kappa<\lambda$, let $g:\theta\rightarrow\lambda$, and let $h\in\prod A^{\bar{c}}$. Fix $\beta<\theta$ and write $g^{\bar{c}}(\beta)=X_{\alpha\beta}$. Assume that this cell is not small. If it is large, assume in addition that
  $$\hat{g}(\kappa_{\alpha\beta})<h(\kappa_{\alpha\beta}).$$
  Then
  $$|g^{\bar{c}\downarrow h}(\beta)|<|g^{\bar{c}}(\beta)|.$$
\end{lemma}
\begin{proof}
If $X_{\alpha\beta}$ is medium, the conclusion follows from the definition of $h$-refinement. If it is large, the definition of $\hat{g}$ and the assumed bound ensure that $g(\beta)\in X_{\alpha\beta}[h(\kappa_{\alpha\beta})]$, whose cardinality is smaller than $|X_{\alpha\beta}|$.
\end{proof}

We finish with an elementary closure observation that will be used without further comment in the later proofs.

\begin{observation}
  \label{partitionclosure}
  Suppose $\theta\leq\kappa<\lambda$ are infinite cardinals and $\mu$ is an infinite cardinal. For each $(\kappa,\theta)$-partition sequence $\bar c$ of $\lambda$, let $\mathcal H_{\bar c}\subseteq\prod A^{\bar c}$ have cardinality at most $\mu$. Then there is a family $\mathcal S$ of such partition sequences, of cardinality at most $\mu$, containing the trivial partition sequence and satisfying
  $$\bar c\in\mathcal S,\ h\in\mathcal H_{\bar c}\Longrightarrow\bar c\downarrow h\in\mathcal S.$$
\end{observation}
\begin{proof}
  Let $\mathcal S_0$ consist of the trivial partition sequence, in which each $c_\beta$ is constantly zero, and recursively define
  $$\mathcal S_{n+1}=\mathcal S_n\cup\{\bar c\downarrow h:\bar c\in\mathcal S_n,\ h\in\mathcal H_{\bar c}\}.$$
  Each $\mathcal S_n$ has cardinality at most $\mu$, so $\mathcal S=\bigcup_{n<\omega}\mathcal S_n$ has the required closure property and cardinality at most $\mu$.
\end{proof}

In each of the following sections, we will supply the refining families and use the resulting closed family of partition sequences to obtain the required covering.

\section{Characterizing covering numbers}
\label{sec:coveringnumbers}

We now combine the partition machinery of Section~\ref{sec:partitions} with pcf theory to prove Shelah's $\cov$ vs. $\pp$ theorem and characterize the covering number in terms of small products, pseudopowers, and the full product $\prod \mathsf{R}(\kappa,\lambda]$. Our proof bypasses the Skolem hulls and characteristic functions used by previous proofs of this theorem.

We start with a reformulation of Claim~3.2 of Chapter~IX of~\cite{cardarith}. This lemma translates the local product cofinalities needed below into statements about standard pcf theory. Analogous local characterizations will be used in later sections, and all will have similar proofs grounded in the pcf theory of products of progressive sets of cardinals.  Note that this lemma also applies when $\sigma=\aleph_0$, although other parts of our proof will require the uncountability assumption. 

\begin{lemma}
  \label{conventional}
  Suppose $A$ is a progressive set of regular cardinals and $\sigma=\cf(\sigma)<\theta$. Then
  $$\cf_{\Gamma(\theta,\sigma)}\left(\prod A\right)=\sup\pcf_{\Gamma(\theta,\sigma)}(A).$$
\end{lemma}
\begin{proof}
  Let $\mu=\sup\pcf_{\Gamma(\theta,\sigma)}(A).$  Suppose $\mathscr{F}\subseteq\prod A$ is $\Gamma(\theta,\sigma)$-cofinal, and let
  $$\tau=\tcf(\prod B/J)$$
  witness a member of $\pcf_{\Gamma(\theta,\sigma)}(A)$. If $|\mathscr{F}|<\tau$, then the restrictions of the members of $\mathscr{F}$ to $B$ have a common upper bound $g\in\prod B$ modulo $J$. Choose $\mathscr{F}_0\in[\mathscr{F}]^{<\sigma}$ such that for every $b\in B$ there is an $f\in\mathscr{F}_0$ with $g(b)<f(b)$. For each $f\in\mathscr{F}_0$, the set
  $$\{b\in B:g(b)<f(b)\}$$
  belongs to $J$, and their union is $B$, contradicting the $\sigma$-completeness and properness of $J$. Thus, every $\Gamma(\theta,\sigma)$-cofinal family has cardinality at least $\mu$.

  Conversely, fix a pcf generating sequence for $A$. For each $\tau\in\pcf_{\Gamma(\theta,\sigma)}(A)$, choose a family $\mathscr{F}_\tau\subseteq\prod A$ of cardinality $\tau$ whose restrictions to $B_\tau[A]$ are pointwise cofinal in $\prod B_\tau[A]$, and let
  $$\mathscr{F}=\bigcup\{\mathscr{F}_\tau:\tau\in\pcf_{\Gamma(\theta,\sigma)}(A)\}.$$
  Then $|\mathscr{F}|\leq\mu$. If $\mathscr{F}$ fails to be $\Gamma(\theta,\sigma)$-cofinal, choose $B\in[A]^{<\theta}$ and $g\in\prod B$ witnessing this. Let $\hat J$ be the $\sigma$-complete ideal on $A$ generated by $A\setminus B$ together with the sets
  $$X_f=\{b\in B:g(b)<f(b)\}$$
  for $f\in\mathscr{F}$. No union of fewer than $\sigma$ of the sets $X_f$ covers $B$, so $\hat J$ is proper.

  Let $\tau\in\pcf(A)$ be least such that $B_\tau[A]\notin\hat J$, and define
  $$X\in K\Longleftrightarrow X\cap B_\tau[A]\in\hat J.$$
  Then $K$ is a proper $\sigma$-complete ideal extending $\hat J$ and containing $A\setminus B_\tau[A]$. Since $J_{<\tau}[A]\subseteq\hat J$, the generator theorem gives
  $$\tau=\tcf(\prod A/K).$$
  Since $A\setminus B\in K$, this is equivalently a representation over $B$, and therefore $\tau\in\pcf_{\Gamma(\theta,\sigma)}(A)$. Extend $g$ to $\tilde g\in\prod A$. Applying pointwise cofinality to $\tilde g+1$, choose $f\in\mathscr{F}_\tau$ such that $\tilde g(a)<f(a)$ for every $a\in B_\tau[A]$. Then
  $$B_\tau[A]\subseteq(A\setminus B)\cup X_f\in\hat J,$$
  a contradiction.
\end{proof}

The next corollary applies Lemma~\ref{partitioncriterion} to $\Gamma(\theta^+,\sigma)$-cofinal families.

\begin{corollary}
  \label{partitionreduction}
  Suppose $\bar{c}$ is a $(\kappa,\theta)$-partition sequence of $\lambda$, where $\theta\leq\kappa<\lambda$, and suppose $\sigma\leq\theta$ is regular. If $\mathcal{H}$ is a nonempty $\Gamma(\theta^+,\sigma)$-cofinal subset of $\prod A^{\bar{c}}$, then for every $g:\theta\rightarrow\lambda$ there is an $\mathcal{H}_0\in[\mathcal{H}]^{<\sigma}$ such that $\{\bar{c}\downarrow h:h\in\mathcal{H}_0\}$ reduces $g$ over $\bar{c}$.
\end{corollary}
\begin{proof}
  Fix $g:\theta\rightarrow\lambda$, and let
  $$B=\{\kappa_{\alpha\beta}:g^{\bar{c}}(\beta)=X_{\alpha\beta}\text{ is large in }\bar{c}\}.$$
  Since $|B|\leq\theta$, apply $\Gamma(\theta^+,\sigma)$-cofinality to $\hat{g}\upharpoonright B$ to obtain $\mathcal H_0\in[\mathcal H]^{<\sigma}$. If this family is empty, enlarge it by one member of $\mathcal H$. At each coordinate where the selected cell is not small, Lemma~\ref{partitioncriterion} gives a shrinking refinement: any member of $\mathcal H_0$ works for a medium cell, and cofinality supplies a suitable member for a large cell. Thus $\{\bar c\downarrow h:h\in\mathcal H_0\}$ reduces $g$ over $\bar c$.
\end{proof}

The next proposition applies Lemma~\ref{partitiontree} to ordinary coverings.

\begin{proposition}
  \label{partitioncover}
  Suppose $\aleph_0<\sigma\leq\theta\leq\kappa<\lambda$, with $\sigma$ regular. If $\mathcal{S}$ is a nonempty collection of $(\kappa,\theta)$-partition sequences of $\lambda$ such that for every $\bar{c}\in\mathcal{S}$ and every $g:\theta\rightarrow\lambda$ there is an $\mathcal{S}_0\in[\mathcal{S}]^{<\sigma}$ that reduces $g$ over $\bar{c}$, then
  $$\cov(\lambda,\kappa^+,\theta^+,\sigma)\leq|\mathcal{S}|.$$
\end{proposition}
\begin{proof}
  For each $\bar{c}\in\mathcal{S}$, let $C^{\bar{c}}$ be the union of its small cells, and set
  $$\mathcal{P}=\{C^{\bar{c}}:\bar{c}\in\mathcal{S}\}.$$
  Then $\mathcal{P}\subseteq[\lambda]^{\leq\kappa}$ and $|\mathcal{P}|\leq|\mathcal{S}|$.

  Fix $g:\theta\rightarrow\lambda$. Starting with any member of $\mathcal{S}$ at the root, recursively construct a tree $T\subseteq{}^{<\omega}\sigma$ and label its nodes by partition sequences $\bar{c}^s\in\mathcal{S}$. At a node $s$, choose a family $\mathcal{S}_s\in[\mathcal{S}]^{<\sigma}$ that reduces $g$ over $\bar{c}^s$, and label the immediate successors of $s$ by the members of $\mathcal{S}_s$. The assigned partition sequences refine along branches, and the shrinking hypothesis of Lemma~\ref{partitiontree} holds at every coordinate $\beta<\theta$.

  Each node has fewer than $\sigma$ immediate successors. Regularity of $\sigma$ implies that every level of $T$ has cardinality less than $\sigma$, and uncountability of $\sigma$ then gives $|T|<\sigma$. Applying Lemma~\ref{partitiontree} at each coordinate, we obtain
  $$\ran(g)\subseteq\bigcup_{s\in T}C^{\bar{c}^s}.$$
  Thus, fewer than $\sigma$ members of $\mathcal{P}$ cover the range of $g$. Every set in $[\lambda]^{\leq\theta}$ is contained in the range of some such function, so $\mathcal{P}$ witnesses the desired inequality.
\end{proof}

The next theorem is our treatment of Shelah's $\cov$ vs. $\pp$ theorem.   We have limited ourselves to a particular configuration of the parameters, but the full theorem (as presented in Chapter~II of~\cite{cardarith}) can be derived from the following result fairly easily.  Our proof gives the added equivalences of (2) and (4), showing that under our hypotheses, the same phenomenon registers in four different contexts,  

\begin{theorem}[$\cov$ vs. $\pp$]
  \label{covpp}
  Suppose $\aleph_0<\sigma\leq\theta\leq\kappa<\lambda$ with $\sigma$ regular. Then the following cardinals are equal:
  \begin{enumerate}
    \item $\cov(\lambda,\kappa^+,\theta^+,\sigma)$.

      \sk

    \item $\sup\{\sup\pcf_{\Gamma(\theta^+,\sigma)}(A):A\subseteq \mathsf{R}(\kappa,\lambda]\wedge |A|\leq\kappa\}$.

      \sk

    \item $\lambda+\sup\{\pp_{\Gamma(\theta^+,\sigma)}(\nu): \sigma\leq\cf(\nu)\leq\theta\leq\kappa<\nu\leq\lambda\}$.

      \sk

    \item $\cf_{\Gamma(\theta^+,\sigma)}\bigl(\prod \mathsf{R}(\kappa,\lambda]\bigr)$.

  \end{enumerate}
\end{theorem}

\begin{proof}
  Let $R=\mathsf{R}(\kappa,\lambda]$, and let $\mu_i$ denote the cardinal in item $(i)$. We prove
  $$\mu_1\leq\mu_2\leq\mu_3\leq\mu_4\leq\mu_1.$$
  To prove $\mu_1\leq\mu_2$, note that singleton sets $A$ give $\lambda=\sup R\leq\mu_2$.\footnote{We use the conventions $\pcf_{\Gamma(\theta,\sigma)}(\varnothing)=\varnothing$ and $\cf_{\Gamma(\theta,\sigma)}(\prod\varnothing)=0$.} For each $(\kappa,\theta)$-partition sequence $\bar c$ of $\lambda$ with $A^{\bar c}\ne\varnothing$, Lemma~\ref{conventional} gives a nonempty $\Gamma(\theta^+,\sigma)$-cofinal family
  \[
  \mathcal H_{\bar c}\subseteq\prod A^{\bar c}
  \]
  of cardinality at most $\mu_2$. If $A^{\bar c}=\varnothing$, take $\mathcal H_{\bar c}=\{\varnothing\}$.

Next, we shoose a family $\mathcal S$ of partition sequences of cardinality at most $\mu_2$, containing the trivial partition sequence and closed under the operation $\bar c\mapsto\bar c\downarrow h$ for $h\in\mathcal H_{\bar c}$.

Given a partition sequence $\bar c\in\mathcal S$ and function $g:\theta\rightarrow\lambda$, we can apply Corollary~\ref{partitionreduction} and conclude there is a  family
  $\mathcal H_0\subseteq\mathcal H_{\bar c}$ of cardinality less than
  $\sigma$ such that the collection $\{\bar c\downarrow h:h\in\mathcal H_0\}$  reduces $g$ over $\bar c$. Every member of this family belongs
  to $\mathcal S$ by closure. Thus $\mathcal S$ satisfies
  the hypothesis of Proposition~\ref{partitioncover}, and consequently
  \[
  \mu_1=\operatorname{cov}(\lambda,\kappa^+,\theta^+,\sigma)
  \leq|\mathcal S|
  \leq\mu_2.
  \]

To prove $\mu_2\leq\mu_3$, we use a standard argument from pcf-theory that we will spell out carefully for the benefit of the reader.  Consider a representation $\tau=\tcf(\prod B/J)$ contributing to $\mu_2$, where $B\subseteq R$, $|B|\leq\theta$, and $J$ is a proper $\sigma$-complete ideal on $B$. If $\tau\leq\lambda$, then $\tau\leq\mu_3$. Assume, therefore, that $\tau>\lambda$.  It follows that $J$ must contain every singleton from $B$, and $B\cap\lambda$ must be $J$-positive (even in the dual filter!).  We can let $\nu\leq\lambda$ be the least ordinal such that $C:=B\cap\nu$ is $J$-positive. Since $J$ must contains  singletons $\nu$ cannot be a successor ordinal.  By minimality, $C$ is cofinal in $\nu$, and every bounded subset of $C$ belongs to $J\upharpoonright C$. Since $C$ consists of cardinals and $\nu$ is a limit ordinal, $\nu$ is a limit cardinal. Also, $C$ is nonempty and its members are greater than $\kappa$, so $\kappa<\nu$.

If $\cf(\nu)<\sigma$, then $C$ would be the union of fewer than $\sigma$ bounded sets, contrary to the $\sigma$-completeness and properness of $J\upharpoonright C$. Consequently,
  $$\sigma\leq\cf(\nu)\leq|C|\leq\theta\leq\kappa<\nu,$$
  so $\nu$ is singular.  Restricting a scale for $\prod B/J$ to the $J$-positive set $C$ preserves both strict increase and cofinality, and therefore
  $$\tau=\tcf\left(\prod C/(J\upharpoonright C)\right).$$
Thus, $\tau$ contributes to $\pp_{\Gamma(\theta^+,\sigma)}(\nu)$, and hence $\tau\leq\mu_3$. Taking suprema gives $\mu_2\leq\mu_3$.

The inequality $\mu_3\leq\mu_4$ is also very easy. Let $\mathscr{F}\subseteq\prod R$ be $\Gamma(\theta^+,\sigma)$-cofinal in $\prod R$ of cardinality $\mu_4$. Singleton coordinates give
  $$\lambda=\sup R\leq\mu_4.$$
 Now suppose $\kappa<\nu\leq\lambda$ with $\sigma\leq\cf(\nu)\leq\theta$, and let
  $$\tau=\tcf(\prod B/J)$$
contribute to $\pp_{\Gamma(\theta^+,\sigma)}(\nu)$. Discarding the bounded initial segment $B\cap(\kappa+1)$ preserves true cofinality, so we may assume $B\subseteq R$ and clearly our choice of $B$ guarantees $|B|\leq\theta$.

If $\mu_4<\tau$, the restrictions of the members of $\mathscr{F}$ to $B$ have a common upper bound $g\in\prod B$ modulo $J$. By $\Gamma(\theta^+,\sigma)$-cofinality, there is an $\mathscr{F}_0\in[\mathscr{F}]^{<\sigma}$ such that
  $$B=\bigcup_{f\in\mathscr{F}_0}\{b\in B:g(b)<f(b)\}.$$
Each comparison set belongs to $J$, contradicting the $\sigma$-completeness and properness of $J$. Thus every such $\tau$ is at most $\mu_4$, and hence $\mu_3\leq\mu_4$.

  To prove $\mu_4\leq\mu_1$, let $\mathcal{P}\subseteq[\lambda]^{\leq\kappa}$ witness $\mu_1$. For each $X\in\mathcal{P}$ define
  $$f_X(\rho)=\sup\bigl(\{\alpha+1:\alpha\in X\cap\rho\}\cup\{0\}\bigr)$$
  for $\rho\in R$. Since $|X|\leq\kappa<\rho$, we have $f_X(\rho)<\rho$. We claim that
  $$\mathscr{F}=\{f_X:X\in\mathcal{P}\}$$
  is $\Gamma(\theta^+,\sigma)$-cofinal. Indeed, suppose
  $$B\in[R]^{<\theta^+}$$
  and $g\in\prod B$. The range of $g$ has cardinality at most $\theta$, so there is a family $\mathcal{P}_0\in[\mathcal{P}]^{<\sigma}$ covering it. For each $\rho\in B$, choose $X\in\mathcal{P}_0$ with $g(\rho)\in X$; then
  $$g(\rho)<f_X(\rho).$$
  Thus, $\mu_4\leq\mu_1$, and all four cardinals are therefore equal.
\end{proof}

\section{Coverings and normal ideals}
\label{sec:normalideals}

We now consider the corresponding covering problem for normal ideals.  In this setting the role of a union of fewer than $\sigma$ sets in the preceding section will be taken by a diagonal union. Recall that if $\kappa$ is an uncountable regular cardinal and $\langle X_i:i<\kappa\rangle$ is a sequence of subsets of $\kappa$, its diagonal union is
$$\nabla_{i<\kappa}X_i=\{\delta<\kappa:(\exists i<\delta)\ \delta\in X_i\}.$$
An ideal on $\kappa$ is \emph{normal} if it contains the bounded subsets of $\kappa$ and is closed under diagonal unions. 

Some of the definitions in this section trace back to Section 4 of Chapter~VIII of~\cite{cardarith} and the first part of~\cite{moreca}.  There, Shelah is concerned with a more abstract covering property he calls $\prc$, but here we restrict our attention to the concrete case of measuring cofinality modulo normal ideals.  We develop the material locally in two steps: first for sequences of regular cardinals and normal ideals on $\kappa$, and then for progressive sets of regular cardinals without a fixed enumeration.

\begin{definition}
Suppose $\kappa$ is an uncountable regular cardinal and let $\bar{\kappa}=\langle\kappa_i:i<\kappa\rangle$ be a progressive sequence of regular cardinals, that is, $\kappa<\min\{\kappa_i:i<\kappa\}$. A set of functions $\mathscr{F}\subseteq\prod\bar{\kappa}$ is {\em diagonally cofinal } in $\prod\bar{\kappa}$ if for every $h\in\prod \bar{\kappa}$ there is a sequence $\langle f_\alpha:\alpha<\kappa\rangle$ of functions from $\mathscr{F}$ such that
  $$\{\delta<\kappa: \sup_{\alpha<\delta}f_\alpha(\delta)\leq h(\delta)\}\text{ is nonstationary.}$$
\end{definition}

Diagonally cofinal families in $\prod\bar{\kappa}$ are precisely those collections of functions that remain unbounded modulo any normal ideal on $\kappa$:

\begin{proposition}
  \label{propnorcfchar}
  Suppose $\kappa$ is an uncountable regular cardinal and $\bar{\kappa}$ is a sequence $\langle\kappa_i:i<\kappa\rangle$ of regular cardinals.  Then the following are equivalent for a family $\mathscr{F}\subseteq\prod\bar{\kappa}$:
  \begin{enumerate}
      \sk
    \item $\mathscr{F}$ is diagonally cofinal in $\prod\bar{\kappa}$.
      \sk
    \item $\mathscr{F}$ is unbounded in $\prod\bar{\kappa}/J$ for every proper normal ideal $J$ on $\kappa$.
  \end{enumerate}
\end{proposition}
\begin{proof}
  Suppose first that $\mathscr{F}$ is diagonally cofinal, and let $J$ be a proper normal ideal on $\kappa$. To see that $\mathscr{F}$ is unbounded modulo $J$, fix $h\in\prod\bar{\kappa}$. Choose a sequence $\langle f_\alpha:\alpha<\kappa\rangle$ from $\mathscr{F}$ and a club $C\subseteq\kappa$ such that
  $$\delta\in C\Longrightarrow h(\delta)<\sup_{\alpha<\delta}f_\alpha(\delta).$$
  For each $\alpha<\kappa$, let
  $$X_\alpha=\{\delta<\kappa:h(\delta)<f_\alpha(\delta)\}.$$
  Then
  $$C\subseteq\nabla_{\alpha<\kappa}X_\alpha.$$
  If every $X_\alpha$ belonged to $J$, normality would imply that the diagonal union on the right belonged to $J$. This is impossible, since a proper normal ideal contains every nonstationary set and hence no club. Thus, $X_\alpha$ is $J$-positive for some $\alpha<\kappa$, so $f_\alpha\nleq_J h$. Since $h$ was arbitrary, $\mathscr{F}$ is unbounded in $\prod\bar{\kappa}/J$.

  Conversely, suppose $\mathscr{F}$ is not diagonally cofinal. Fix $h\in\prod\bar{\kappa}$ witnessing this, so that for every sequence $\langle f_\alpha:\alpha<\kappa\rangle$ from $\mathscr{F}$, the set
  $$\left\{\delta<\kappa:\sup_{\alpha<\delta}f_\alpha(\delta)\leq h(\delta)\right\}$$
  is stationary. For each $f\in\mathscr{F}$, define
  $$X_f=\{\delta<\kappa:h(\delta)<f(\delta)\},$$
  and let $J$ be the normal ideal generated by $\NS_\kappa$ together with the sets $X_f$ for $f\in\mathscr{F}$.

  We claim that $J$ is proper. Otherwise, there would be a sequence $\langle f_\alpha:\alpha<\kappa\rangle$ from $\mathscr{F}$ and a nonstationary set $N\subseteq\kappa$ such that
  $$\kappa=N\cup\nabla_{\alpha<\kappa}X_{f_\alpha}.$$
  Choose a club $C\subseteq\kappa$ disjoint from $N$. Then for every $\delta\in C$ there is an $\alpha<\delta$ such that
  $$h(\delta)<f_\alpha(\delta),$$
  contradicting the choice of $h$. Therefore, $J$ is a proper normal ideal. But $X_f\in J$ for every $f\in\mathscr{F}$, so $f\leq_J h$ for every $f\in\mathscr{F}$. Thus, $\mathscr{F}$ is bounded in $\prod\bar{\kappa}/J$, completing the proof.
\end{proof}

Now we can make the natural definition associated with diagonally cofinal subsets of $\prod\bar{\kappa}$.  The following is equivalent to that given as Definition 4.4 on page 343 of~\cite{cardarith}, and better captures a sense of what is being measured.

\begin{definition}
  Suppose $\kappa$ and $\bar{\kappa}$ are as above. Then $\operatorname{nor-cf}\left(\prod\bar\kappa\right)$ is the minimum cardinality of a diagonally cofinal subset of $\prod\bar\kappa$.
\end{definition}

In general, $\norcf(\prod\bar\kappa)$ is going to depend on the enumeration $\bar{\kappa}$, although it will be unchanged if the sequence is modified only on a nonstationary set.  We can, however, still characterize it in terms of the pcf structure of the product $\prod_{i<\kappa}\kappa_i$.

\begin{definition}
Let $\bar\kappa=\langle\kappa_i:i<\kappa\rangle$ be a progressive sequence of regular cardinals, and let $\langle B_\tau[\bar\kappa]:\tau\in\pcf(\bar\kappa)\rangle$ be a generating sequence.\footnote{In this indexed context, the generators are subsets of the index set $\kappa$, and we choose them so that  $B_\Lambda[\bar\kappa]=\kappa$ when $\Lambda=\max\pcf(\bar\kappa)$.}
\begin{enumerate}
  \item  If $I$ is an ideal on $\kappa$, let $\operatorname{Ncl}(I)$ denote the normal ideal generated by $I$. Thus, $X\in\operatorname{Ncl}(I)$ if and only if there are a nonstationary set $N\subseteq\kappa$ and a sequence $\langle X_\alpha:\alpha<\kappa\rangle$ from $I$ such that

  $$X\subseteq N\cup\nabla_{\alpha<\kappa}X_\alpha.$$
\sk 
\item Finally, define $\pcfnor(\bar\kappa)=\{\tau\in\pcf(\bar\kappa):B_\tau[\bar\kappa]\notin\operatorname{Ncl}(J_{<\tau}[\bar\kappa])\}.$
  \end{enumerate}
\end{definition}

Before moving on, we note that standard pcf arguments tell us that the definition of $\pcfnor$ we give above yields exactly the set of true cofinalities obtained from $\prod\bar{\kappa}$ modulo proper normal ideals on $\kappa$.  Indeed, if
  $$\tau=\tcf\left(\prod\bar\kappa/J\right)$$
  for such an ideal $J$, then the generator theorem gives $J_{<\tau}[\bar\kappa]\subseteq J$ and $B_\tau[\bar\kappa]\notin J$. Since $J$ is normal, it contains $\operatorname{Ncl}(J_{<\tau}[\bar\kappa])$, and hence $\tau\in\pcfnor(\bar\kappa)$.

  Conversely, suppose $\tau\in\pcfnor(\bar\kappa)$, and set
  $$I_\tau=\operatorname{Ncl}(J_{<\tau}[\bar\kappa]).$$
  Define an ideal $J$ on $\kappa$ by
  $$X\in J\Longleftrightarrow X\cap B_\tau[\bar\kappa]\in I_\tau.$$
  Then $J$ is a proper normal ideal extending $J_{<\tau}[\bar\kappa]$ and containing $\kappa\setminus B_\tau[\bar\kappa]$. The generator theorem gives
  $$\tau=\tcf\left(\prod\bar\kappa/J\right).$$
With this observation,  we are now in a position to characterize $\norcf(\prod\bar{\kappa})$. 
\begin{lemma}
  \label{normalcofinality}
  Suppose $\kappa$ is an uncountable regular cardinal and $\bar\kappa=\langle\kappa_i:i<\kappa\rangle$ is a progressive sequence of regular cardinals. Then 
  \begin{equation}
    \label{eqn1}
      \norcf\left(\prod\bar\kappa\right)=\sup\pcfnor(\bar\kappa)
  \end{equation}
\end{lemma}
\begin{proof}
Let $\mu=\sup\pcfnor(\bar{\kappa})$.  Proposition~\ref{propnorcfchar} gives $\norcf(\prod\bar\kappa)\geq\mu$: every diagonally cofinal family is unbounded modulo each proper normal ideal $J$, while a family of cardinality less than $\tcf(\prod\bar\kappa/J)$ is bounded modulo $J$.

  For the reverse inequality, choose for each $\tau\in\pcfnor(\bar\kappa)$ a family $\mathscr{F}_\tau\subseteq\prod\bar\kappa$ of cardinality $\tau$ whose restrictions to $B_\tau[\bar\kappa]$ are pointwise cofinal in $\prod_{i\in B_\tau[\bar\kappa]}\kappa_i$. Set
  $$\mathscr{F}=\bigcup_{\tau\in\pcfnor(\bar\kappa)}\mathscr{F}_\tau.$$
    There are at most $\mu$ cardinals in $\pcfnor(\bar\kappa)$, and each $\mathscr{F}_\tau$ has cardinality at most $\mu$, so $|\mathscr{F}|\leq\mu$.

By Proposition~\ref{propnorcfchar}, it suffices to show that $\mathscr{F}$ is unbounded modulo every proper normal ideal on $\kappa$.  Fix such an ideal $J$, let $h$ be a function in $\prod\bar{\kappa}$, and let $\tau$ be least such that $B_\tau[\bar{\kappa}]\notin J$.

By the minimality of $\tau$, $J_{<\tau}[\bar\kappa]\subseteq J$, so normality gives
  $$\operatorname{Ncl}(J_{<\tau}[\bar\kappa])\subseteq J,$$
and therefore $\tau\in\pcfnor(\bar\kappa)$. Applying pointwise cofinality to $h+1$, choose $f\in\mathscr{F}_\tau$ with $h(i)<f(i)$ for every $i\in B_\tau[\bar\kappa]$. The set $\{i<\kappa:h(i)<f(i)\}$ is therefore $J$-positive and we are done. 
\end{proof}

We next modify these concepts to remove dependence on a particular enumeration.  The following does the trick, and the natural characterization follows by arguments similar to what we saw above.

\begin{definition}
Suppose $\kappa$ is an uncountable regular cardinal and $A$ is a set of regular cardinals with $\kappa<\min(A)$. 
\begin{enumerate}
  \item A family $\mathscr{F}\subseteq\prod A$ is \emph{normally $\kappa$-cofinal} if, for every function $\bar\kappa:\kappa\rightarrow A$, the family $\{f\circ\bar\kappa:f\in\mathscr{F}\}$
  is diagonally cofinal in $\prod\bar{\kappa}$.
  \sk 
  \item A cardinal $\tau$ is in $\pcfnor(A)$ if there is a function  $\bar{\kappa}:\kappa\rightarrow A$ such that $\tau\in\pcfnor(\bar{\kappa})$.
\end{enumerate}
\end{definition}

\begin{lemma}
  \label{normalimages}
  Suppose $\kappa$ is an uncountable regular cardinal and $A$ is a set of regular cardinals with $|A|\leq\kappa<\min(A)$. Then the following two cardinals are equal:
  \begin{enumerate}
    \item the minimum cardinality of a normally $\kappa$-cofinal family in $\prod A$;
      \sk
    \item $\sup\pcfnor(A)$.
  \end{enumerate}
\end{lemma}
\begin{proof}
  The proof is a routine modification of that of Lemma~\ref{normalcofinality}, using the ideals on $A$ induced by maps $\kappa\rightarrow A$, and is left to the reader.
\end{proof}

Moving on to pseudopowers, we arrive at a notion that has been studied in several contexts by Shelah.  

\begin{definition}
  Let $\nu$ be singular with $\kappa=\cf(\nu)>\aleph_0$. Then $\ppnor(\nu)$ is the supremum of all cardinals of the form $\tcf(\prod_{i<\kappa}\nu_i/J)$, where $\langle\nu_i:i<\kappa\rangle$ is a sequence of regular cardinals from $(\kappa,\nu)$, $J$ is a proper normal ideal on $\kappa$, and the sequence converges to $\nu$ modulo $J$, that is,
  $$\{i<\kappa:\nu_i\leq\alpha\}\in J$$
  for every $\alpha<\nu$.
\end{definition}

We would like to use the above pcf-theoretic information as part of a cov vs. pp theorem that focuses on normal ideals.  In order for that to work, we need to think about what sort of covering number might correspond to $\pcfnor$ and $\ppnor$, and we are led to the following definition: 

\begin{definition}
  A family $\mathscr{P}\subseteq[\lambda]^\kappa$ is \emph{diagonally cofinal} if, for every $g:\kappa\rightarrow\lambda$, there is a sequence $\langle X_\alpha:\alpha<\kappa\rangle$ from $\mathscr{P}$ such that
  $$\left\{\delta<\kappa:g(\delta)\notin\bigcup_{\alpha<\delta}X_\alpha\right\}$$
  is nonstationary.
\end{definition}

We are now in a position to state and prove an analogue of Theorem~\ref{covpp} tailored to normal ideals.  The following theorem is a relative of Lemma~1.1 of~\cite{moreca}, where Shelah characterizes a cardinal he calls $\prc_J$. The argument given there is quite sketchy, but it does provide the seed of the idea we are developing in this paper. 

\begin{theorem}
  \label{norcovpp}
  Suppose $\kappa<\lambda$ with $\kappa$ an uncountable regular cardinal.  Then the following four cardinals are equal:
  \begin{enumerate}
    \item The minimum cardinality of a diagonally cofinal subset $\mathscr{P}$ of $[\lambda]^\kappa$.
      \sk
    \item The supremum of all cardinals $\sup\pcfnor(A)$, where $A$ is a nonempty subset of $\mathsf{R}(\kappa,\lambda]$ of cardinality at most $\kappa$.
      \sk
    \item $\lambda+\sup\{\ppnor(\nu):\kappa<\nu\leq\lambda\wedge\cf(\nu)=\kappa\}$.
      \sk
    \item The minimum cardinality of a normally $\kappa$-cofinal family in $\prod \mathsf{R}(\kappa,\lambda]$.
  \end{enumerate}
\end{theorem}
\begin{proof}
  Let $R=\mathsf{R}(\kappa,\lambda]$, and let $\mu_i$ denote the cardinal in item $(i)$. We prove
  $$\mu_1\leq\mu_2\leq\mu_3\leq\mu_4\leq\mu_1.$$

  To prove $\mu_1\leq\mu_2$, first note that constant sequences give $\rho\in\pcfnor(\{\rho\})$ for every $\rho\in R$, and hence
  $$\lambda=\sup R\leq\mu_2.$$
  For each $(\kappa,\kappa)$-partition sequence $\bar c$ of $\lambda$ with $A^{\bar c}\ne\varnothing$, use Lemma~\ref{normalimages} to choose a normally $\kappa$-cofinal family $\mathscr{F}_{\bar c}\subseteq\prod A^{\bar c}$ of cardinality at most $\mu_2$. If $A^{\bar c}=\varnothing$, set $\mathscr{F}_{\bar c}=\{\varnothing\}$. Choose a family $\mathcal{S}$ of cardinality at most $\mu_2$ containing the trivial partition sequence and closed under $\bar c\mapsto\bar c\downarrow f$ for $f\in\mathscr{F}_{\bar c}$.

  Let $C^{\bar c}$ be the union of the small cells of $\bar c$, and put
  $$\mathscr{P}=\{C^{\bar c}\cup\kappa:\bar c\in\mathcal{S}\}\subseteq[\lambda]^\kappa.$$
  Then $|\mathscr{P}|\leq\mu_2$. We show that $\mathscr{P}$ is diagonally cofinal.

  Fix $g:\kappa\rightarrow\lambda$. Starting with the trivial partition sequence at the root, construct partition sequences $\bar c^s\in\mathcal{S}$ for $s\in{}^{<\omega}\kappa$. At a node $s$, set $A^s=A^{\bar c^s}$. If $A^s$ is nonempty, let $\hat g^s\in\prod A^s$ be the function from Definition~\ref{partitionbound}. Set $\eta^s_\beta=\cf(|g^{\bar c^s}(\beta)|)$ when this cell is large, and $\eta^s_\beta=\min(A^s)$ otherwise. Normal $\kappa$-cofinality of $\mathscr{F}_{\bar c^s}$ gives a sequence $\langle f^s_\xi:\xi<\kappa\rangle$ from this family and a club $D_s\subseteq\kappa$ such that
  $$\delta\in D_s\Longrightarrow\hat g^s(\eta^s_\delta)<\sup_{\xi<\delta}f^s_\xi(\eta^s_\delta).$$
  If $A^s$ is empty, let every $f^s_\xi$ be the empty function and put $D_s=\kappa$. In either case, set
  $$\bar c^{s^\smallfrown\langle\xi\rangle}=\bar c^s\downarrow f^s_\xi.$$
  These successors belong to $\mathcal{S}$ and refine $\bar c^s$. Put $Y^s=C^{\bar c^s}\cup\kappa\in\mathscr{P}$.

  Fix a bijection $\pi:{}^{<\omega}\kappa\rightarrow\kappa$, and choose a club $E\subseteq\kappa$ of nonzero limit ordinals such that
  $$\pi``({}^{<\omega}\delta)\subseteq\delta$$
  for every $\delta\in E$. Then
  $$D=E\cap\mathop{\triangle}_{\zeta<\kappa}D_{\pi^{-1}(\zeta)}$$
  is club. Set $Y_\zeta=Y^{\pi^{-1}(\zeta)}$. If $\delta\in D$ and $s\in{}^{<\omega}\delta$, then $\pi(s)<\delta$ and $\delta\in D_s$.

  Fix $\delta\in D$ and consider the tree $T={}^{<\omega}\delta$. At a node $s\in T$, if $g^{\bar c^s}(\delta)$ is medium, choose any $\xi<\delta$. If it is large, $\delta\in D_s$ gives a $\xi<\delta$ such that $\hat g^s(\eta^s_\delta)<f^s_\xi(\eta^s_\delta)$. In either case, Lemma~\ref{partitioncriterion} shows that the chosen successor shrinks the cell containing $g(\delta)$. Applying Lemma~\ref{partitiontree} at coordinate $\delta$, we obtain
  $$g(\delta)\in\bigcup_{s\in{}^{<\omega}\delta}C^{\bar c^s}\subseteq\bigcup_{\zeta<\delta}Y_\zeta.$$
  Thus $\mathscr{P}$ is diagonally cofinal and $\mu_1\leq\mu_2$.

  To prove $\mu_2\leq\mu_3$, consider a representation
  $$\tau=\tcf\left(\prod_{i<\kappa}\rho_i/J\right)>\lambda$$
  contributing to $\mu_2$, where each $\rho_i\in R$ and $J$ is a proper normal ideal on $\kappa$. For each $\rho\in R$, the fiber $\{i<\kappa:\rho_i=\rho\}$ belongs to $J$: otherwise, restricting the scale to that fiber would give $\tau=\rho\leq\lambda$. Consequently, there is a least $\nu\leq\lambda$ such that
  $$S=\{i<\kappa:\rho_i<\nu\}\in J^+.$$
  Minimality implies that $\nu$ is a limit cardinal greater than $\kappa$, that $\{\rho_i:i\in S\}$ is cofinal in $\nu$, and that
  $$\{i\in S:\rho_i\leq\alpha\}\in J$$
  for every $\alpha<\nu$. Since $J$ is $\kappa$-complete and the range on $S$ has cardinality at most $\kappa$, we have $\cf(\nu)=\kappa$.

  The ideal
  $$J_S=\{X\subseteq\kappa:X\cap S\in J\}$$
  is proper and normal, and restricting the scale to $S$ preserves true cofinality $\tau$. Replace each $\rho_i$ for $i\notin S$ by a fixed member of $\{\rho_i:i\in S\}$. The resulting sequence, together with $J_S$, witnesses $\tau\leq\ppnor(\nu)\leq\mu_3$. True cofinalities at most $\lambda$ are already bounded by $\mu_3$, so taking suprema gives $\mu_2\leq\mu_3$.

  To prove $\mu_3\leq\mu_4$, let $\mathscr{F}\subseteq\prod R$ witness $\mu_4$. Constant maps give $\lambda\leq\mu_4$. If $\tau=\tcf(\prod\bar\nu/J)$ contributes to $\ppnor(\nu)$ for $\kappa<\nu\leq\lambda$ with $\cf(\nu)=\kappa$, then $\bar\nu:\kappa\rightarrow R$. The pullback of $\mathscr{F}$ along $\bar\nu$ is diagonally cofinal and hence unbounded modulo $J$ by Proposition~\ref{propnorcfchar}. A family of cardinality less than $\tau$ is bounded modulo $J$, so $\tau\leq|\mathscr{F}|=\mu_4$. Taking suprema gives $\mu_3\leq\mu_4$.

  To prove $\mu_4\leq\mu_1$, let $\mathscr{P}\subseteq[\lambda]^\kappa$ witness $\mu_1$. For $X\in\mathscr{P}$ and $\rho\in R$, define
  $$f_X(\rho)=\sup\bigl(\{\gamma+1:\gamma\in X\cap\rho\}\cup\{0\}\bigr).$$
  Since $|X|=\kappa<\rho=\cf(\rho)$, we have $f_X\in\prod R$. Put $\mathscr{F}=\{f_X:X\in\mathscr{P}\}$. Given $\bar\eta:\kappa\rightarrow R$ and $h\in\prod\bar\eta$, diagonal cofinality of $\mathscr{P}$ gives a sequence $\langle X_\alpha:\alpha<\kappa\rangle$ from $\mathscr{P}$ and a club $C\subseteq\kappa$ such that
  $$\delta\in C\Longrightarrow h(\delta)\in\bigcup_{\alpha<\delta}X_\alpha.$$
  For each $\delta\in C$, some $\alpha<\delta$ therefore satisfies
  $$h(\delta)<f_{X_\alpha}(\bar\eta(\delta)).$$
  Thus $\mathscr{F}$ is normally $\kappa$-cofinal and $\mu_4\leq|\mathscr{F}|\leq\mu_1$, completing the cycle.

\end{proof}

\section{Coverings modulo an ideal}
\label{sec:almosteverywhere}

Suppose $\theta\leq\kappa<\lambda$ and $I$ is a proper ideal on $\theta$. For our final example, we now ask how many sets from $[\lambda]^{\leq\kappa}$ are needed so that, for every $g:\theta\rightarrow\lambda$, fewer than $\sigma$ of them cover $g(\beta)$ for $I$-almost every $\beta<\theta$. Our partition method gives a complete answer whenever $\aleph_0<\sigma=\cf(\sigma)\leq\theta$ and $I$ is $\sigma$-complete. In particular, the case $\sigma=\aleph_1$ allows countably many covering sets and requires only countable completeness of $I$. The following definitions describe the corresponding pcf information.

\begin{definition}
  Let $A$ be a set of regular cardinals and let $I$ be a proper ideal on a set $D$ with $|D|<\min(A)$.
  \begin{enumerate}
    \item An ideal $J$ on $A$ is an {\em image of $I$} if there is a sequence $\bar{\eta}=\langle \eta_d:d\in D\rangle$ of cardinals from $A$ such that
      \begin{equation*}
        B\in J\Longleftrightarrow \{d\in D: \eta_d\in B\}\in I.
      \end{equation*}
      \sk
    \item We say $J$ {\em extends an image of $I$} if for some such sequence $\bar{\eta}$ we have
      \begin{equation*}
        B\in J^+\Longrightarrow \{d\in D: \eta_d\in B\}\in I^+.
      \end{equation*}
  \end{enumerate}
\end{definition}

The sequence $\bar\eta$ need not be one-to-one. In Rudin--Keisler terminology, the images of $I$ are exactly the ideals $J\leq_{RK}I$.

We allow extensions because the construction must sometimes restrict attention to an $I$-positive set of coordinates. If $X\in I^+$, then
$$I_X=\{Y\subseteq D:Y\cap X\in I\}$$
is a proper ideal extending $I$, and it is $\sigma$-complete whenever $I$ is $\sigma$-complete. The image of $I_X$ under a sequence $\bar\eta$ extends the image of $I$ under the same sequence. This is the localization used in the covering argument below.

\begin{definition}
  Let $A$, $D$, and $I$ be as above with $A$ progressive, and let $\sigma$ be an infinite regular cardinal. A cardinal $\tau$ is in $\pcf_{\Gamma(I,\sigma)}(A)$ if $\tau=\tcf(\prod A/J)$ for some proper $\sigma$-complete ideal $J$ on $A$ extending an image of $I$, with the true cofinality defined.
\end{definition}

If $J$ is a proper ideal on $A$, we say that $F\subseteq\prod A$ is $<\sigma$-cofinal in $\prod A/J$ if for every $g\in\prod A$ there is an $F_0\in[F]^{<\sigma}$ such that
$$\{\eta\in A:(\forall f\in F_0)\ f(\eta)\leq g(\eta)\}\in J.$$

\begin{definition}
  Let $A$ be a set of regular cardinals, let $I$ be a proper ideal on a set $D$ with $|D|<\min(A)$, and let $\sigma$ be an infinite regular cardinal. A family $F\subseteq\prod A$ is \emph{$\Gamma(I,\sigma)$-cofinal} if $F$ is $<\sigma$-cofinal in $\prod A/J$ whenever $J$ is a proper $\sigma$-complete ideal on $A$ extending an image of $I$. The minimum cardinality of such a family is denoted by $\cf_{\Gamma(I,\sigma)}(\prod A)$.
\end{definition}

\begin{lemma}
  \label{idealcofinality}
  Let $A$ be a progressive set of regular cardinals, let $\sigma$ be an infinite regular cardinal, and let $I$ be a proper ideal on a set $D$ with $|D|<\min(A)$. Then the following two cardinals are equal:
  \begin{enumerate}
    \item $\sup\pcf_{\Gamma(I,\sigma)}(A)$

      \sk

    \item $\cf_{\Gamma(I,\sigma)}(\prod A)$.
  \end{enumerate}
\end{lemma}
\begin{proof}
  The proof is a routine modification of that of Lemma~\ref{conventional}. The class of proper $\sigma$-complete ideals extending an image of $I$ is closed under proper $\sigma$-complete extensions, so the ideals constructed in that argument remain admissible.
\end{proof}

As in Section~\ref{sec:coveringnumbers}, we will also characterize the covering number using pseudopowers. Here the representing ideals must extend an image of $I$.

\begin{definition}
  Let $I$ be a proper ideal on $\theta$, let $\sigma$ be an infinite regular cardinal, and let $\nu$ be singular with $\sigma\leq\cf(\nu)\leq\theta<\nu$. Then $\pp_{\Gamma(I,\sigma)}(\nu)$ is the supremum of all cardinals of the form $\tcf(\prod A/J)$, where $A\subseteq\nu\cap\Reg$ is cofinal in $\nu$, $|A|\leq\theta<\min(A)$, and $J$ is a proper $\sigma$-complete ideal on $A$ extending an image of $I$ and containing every bounded subset of $A$.
\end{definition}

\begin{theorem}
  Suppose $\aleph_0<\sigma\leq\theta\leq\kappa<\lambda$ with $\sigma$ regular, and let $I$ be a proper $\sigma$-complete ideal on $\theta$. Then the following four cardinals are equal:
  \begin{enumerate}
      \sk

    \item The minimum cardinality of a family $\mathcal{C}\subseteq[\lambda]^{\leq\kappa}$ such that, for every $g:\theta\rightarrow\lambda$, there is an $\mathcal{X}\in[\mathcal{C}]^{<\sigma}$ with
      $$\{\beta<\theta:g(\beta)\notin\bigcup\mathcal{X}\}\in I.$$

      \bigskip

    \item The supremum of all cardinals $\sup\pcf_{\Gamma(I,\sigma)}(A)$ where $A$ is a nonempty subset of $\mathsf{R}(\kappa,\lambda]$ of cardinality at most $\kappa$.

      \sk

    \item $\lambda+\sup\{\pp_{\Gamma(I,\sigma)}(\nu):\sigma\leq\cf(\nu)\leq\theta\leq\kappa<\nu\leq\lambda\}$.

      \sk

    \item $\cf_{\Gamma(I,\sigma)}\bigl(\prod \mathsf{R}(\kappa,\lambda]\bigr)$.
  \end{enumerate}
\end{theorem}
\begin{proof}
  Let $R=\mathsf{R}(\kappa,\lambda]$, and let $\mu_i$ denote the cardinal in item $(i)$. We prove
  $$\mu_1\leq\mu_2\leq\mu_3\leq\mu_4\leq\mu_1.$$

  To prove $\mu_1\leq\mu_2$, note that for each $\rho\in R$ the constant map $\theta\rightarrow\{\rho\}$ sends $I$ to the ideal $\{\varnothing\}$ on $\{\rho\}$. Thus $\rho\in\pcf_{\Gamma(I,\sigma)}(\{\rho\})$, and
  $$\lambda=\sup R\leq\mu_2.$$
  For each $(\kappa,\theta)$-partition sequence $\bar c$ of $\lambda$ with $A^{\bar c}\ne\varnothing$, Lemma~\ref{idealcofinality} gives a $\Gamma(I,\sigma)$-cofinal family $\mathcal H_{\bar c}\subseteq\prod A^{\bar c}$ of cardinality at most $\mu_2$. Adjoin the zero function to each such family, and take $\mathcal H_{\bar c}=\{\varnothing\}$ when $A^{\bar c}=\varnothing$. Choose a family $\mathcal S$ of cardinality at most $\mu_2$ containing the trivial partition sequence and closed under these refinements.

  Let $C^{\bar c}$ be the union of the small cells of $\bar c$, and put
  $$\mathcal C=\{C^{\bar c}:\bar c\in\mathcal S\}\subseteq[\lambda]^{\leq\kappa}.$$
  We show that $\mathcal C$ witnesses $\mu_1\leq\mu_2$.

  Fix $g:\theta\rightarrow\lambda$. Starting with the trivial partition sequence at the root, construct a $<\sigma$-branching tree $T\subseteq{}^{<\omega}\sigma$ whose nodes carry partition sequences $\bar c^s\in\mathcal S$ and exceptional sets $E_s\in I$. At a node $s$, set $A^s=A^{\bar c^s}$, $C^s=C^{\bar c^s}$, and
  $$X_s=\{\beta<\theta:g^{\bar c^s}(\beta)\text{ is large in }\bar c^s\}.$$
  If $X_s\in I$, let $H^s$ consist of the zero function on $A^s$, including the empty function when $A^s$ is empty, and put $E_s=X_s$.

  If $X_s\in I^+$, then $A^s$ is nonempty. Set $\eta^s_\beta=\cf(|g^{\bar c^s}(\beta)|)$ for $\beta\in X_s$ and $\eta^s_\beta=\min(A^s)$ otherwise. Define an ideal $J_s$ on $A^s$ by
  $$B\in J_s\Longleftrightarrow\{\beta\in X_s:\eta^s_\beta\in B\}\in I.$$
  This is the image under $\bar\eta^s$ of the proper $\sigma$-complete ideal $I_{X_s}$, so it is proper and $\sigma$-complete and extends the image of $I$ under $\bar\eta^s$. Let $\hat g^s\in\prod A^s$ be the function from Definition~\ref{partitionbound}. Choose a nonempty $H^s\in[\mathcal H_{\bar c^s}]^{<\sigma}$ such that
  $$D_s=\{\rho\in A^s:(\forall h\in H^s)\ h(\rho)\leq\hat g^s(\rho)\}\in J_s,$$
  and put
  $$E_s=\{\beta\in X_s:\eta^s_\beta\in D_s\}\in I.$$

  In either case, enumerate $H^s$ as $\langle h^s_i:i<\delta_s\rangle$, where $0<\delta_s<\sigma$, and set
  $$\bar c^{s^\smallfrown\langle i\rangle}=\bar c^s\downarrow h^s_i$$
  for $i<\delta_s$. Closure ensures that these successors belong to $\mathcal S$. If $\beta\notin E_s$ and the cell containing $g(\beta)$ is large, the definition of $E_s$ gives an $h\in H^s$ with $\hat g^s(\eta^s_\beta)<h(\eta^s_\beta)$. For a medium cell, any member of $H^s$ is suitable. Lemma~\ref{partitioncriterion} therefore shows that whenever $\beta\notin E_s$ and the cell containing $g(\beta)$ is not small, some successor decreases its cardinality.

  Since $\sigma$ is uncountable and regular, $|T|<\sigma$. The $\sigma$-completeness of $I$ therefore gives
  $$E=\bigcup_{s\in T}E_s\in I.$$
  For each $\beta\notin E$, the shrinking hypothesis holds at every node of $T$. Applying Lemma~\ref{partitiontree} at coordinate $\beta$, we obtain
  $$\beta\notin E\Longrightarrow g(\beta)\in\bigcup_{s\in T}C^s.$$
  Each $C^s$ belongs to $\mathcal C$, and $|T|<\sigma$, so $\mathcal C$ has the required covering property. Hence $\mu_1\leq|\mathcal C|\leq\mu_2$.

  For the pseudopower comparisons, note that if $J$ extends an image of $I$ and $B\in J^+$, then $J\upharpoonright B$ also extends an image of $I$. Indeed, let $\langle\eta_\beta:\beta<\theta\rangle$ witness the image condition, choose $b\in B$, and set $\eta'_\beta=\eta_\beta$ when $\eta_\beta\in B$ and $\eta'_\beta=b$ otherwise. For $C\subseteq B$, if $\{\beta<\theta:\eta'_\beta\in C\}\in I$, then $\{\beta<\theta:\eta_\beta\in C\}\in I$, and hence $C\in J\upharpoonright B$.

  To prove $\mu_2\leq\mu_3$, let $\tau=\tcf(\prod A/J)>\lambda$ contribute to $\mu_2$, and let $\bar\eta=\langle\eta_\beta:\beta<\theta\rangle$ witness that $J$ extends an image of $I$. The complement of $\ran(\bar\eta)$ belongs to $J$. Restricting to this range therefore preserves true cofinality, and we may assume $|A|\leq\theta$. As in Section~\ref{sec:coveringnumbers}, every singleton belongs to $J$, and we choose the least $\nu\leq\lambda$ for which $B=A\cap\nu$ is $J$-positive. The same argument gives
  $$\sigma\leq\cf(\nu)\leq|B|\leq\theta\leq\kappa<\nu,$$
  and $J\upharpoonright B$ contains every bounded subset of $B$. It also extends an image of $I$, so restricting the scale to $B$ shows that $\tau$ contributes to $\pp_{\Gamma(I,\sigma)}(\nu)$. True cofinalities at most $\lambda$ are already bounded by $\mu_3$, and hence $\mu_2\leq\mu_3$.

  To prove $\mu_3\leq\mu_4$, let $F\subseteq\prod R$ witness $\mu_4$. Constant maps give $\lambda\leq\mu_4$. Consider a representation
  $$\tau=\tcf(\prod B/J)$$
  contributing to $\pp_{\Gamma(I,\sigma)}(\nu)$, where $\kappa<\nu\leq\lambda$ and $\sigma\leq\cf(\nu)\leq\theta$. Restricting to $B\setminus(\kappa+1)$ preserves true cofinality and the image condition, so we may assume $B\subseteq R$. Extend $J$ to an ideal $\widehat J$ on $R$ by
  $$X\in\widehat J\Longleftrightarrow X\cap B\in J.$$
  Then $\widehat J$ is proper and $\sigma$-complete, extends an image of $I$, and satisfies
  $$\tau=\tcf(\prod R/\widehat J).$$
  If $|F|<\tau$, its members have a common upper bound modulo $\widehat J$. The $<\sigma$-cofinality of $F$ then contradicts the $\sigma$-completeness and properness of $\widehat J$. Thus $\tau\leq|F|=\mu_4$, and taking suprema gives $\mu_3\leq\mu_4$.

  To prove $\mu_4\leq\mu_1$, let $\mathcal{C}$ witness $\mu_1$. For each $X\in\mathcal{C}$, define
  $$h_X(\rho)=\sup(X\cap\rho)+1$$
  for $\rho\in R$. Since $|X|\leq\kappa<\rho=\cf(\rho)$, we have $h_X\in\prod R$. We show that
  $$F=\{h_X:X\in\mathcal{C}\}$$
  is $\Gamma(I,\sigma)$-cofinal in $\prod R$.

  Fix a proper $\sigma$-complete ideal $J$ on $R$ extending an image of $I$, and choose $\bar\eta=\langle\eta_\beta:\beta<\theta\rangle$ witnessing the extension property. Given $f\in\prod R$, apply the covering property to $\beta\mapsto f(\eta_\beta)$ and choose $\mathcal{X}\in[\mathcal{C}]^{<\sigma}$ such that
  $$E=\{\beta<\theta:f(\eta_\beta)\notin\bigcup\mathcal{X}\}\in I.$$
  Put
  $$B_f=\{\rho\in R:(\forall X\in\mathcal{X})\ h_X(\rho)\leq f(\rho)\}.$$
  If $\eta_\beta\in B_f$, then $f(\eta_\beta)$ belongs to none of the sets in $\mathcal{X}$, so $\beta\in E$. Thus the preimage of $B_f$ under $\bar\eta$ belongs to $I$, and hence $B_f\in J$. This proves that $F$ is $\Gamma(I,\sigma)$-cofinal, and therefore
  $$\mu_4\leq|F|\leq|\mathcal{C}|=\mu_1.$$
  This completes the cycle.

\end{proof}

\section{Conclusion}
\label{sec:conclusion}

We conclude by emphasizing how the results of this paper should be viewed. Our main point is not that we have three different theorems characterizing various covering numbers in terms of pcf data. Instead, we want to demonstrate that there is a general technique that works here and that can easily be adapted to other settings: suitable local bounds in products of regular cardinals can be converted into global coverings by an elementary partition argument. The bounds supply refinements that shrink the relevant cells, and a tree of these refinements yields a covering because cell cardinalities cannot decrease indefinitely.

The three applications demonstrate how to reuse this mechanism. For ordinary coverings, diagonal coverings, and coverings modulo a complete ideal, the task is to identify the appropriate notion of cofinality and organize the resulting refinements to meet the covering requirement. The parallel characterizations by local true cofinalities, pseudopowers, and cofinal families in full products arise from this common method.

Finally, in the spirit of the times we add a disclosure on the use of AI in the preparation of this paper:  the use of AI was limited to assistance with editing and proofreading, and the author is fully responsible for the mathematics.  We wrote the paper using the AI capabilities built into the latex editor at {\tt murphy.ai}.

\end{document}